\documentclass[12pt]{amsart}
\usepackage{amsmath}
\usepackage{amssymb}
\usepackage{stmaryrd}
\usepackage{color}

\usepackage{ifthen,fancyhdr}
\newcommand{\notre}[1]{\ifthenelse{#1=13}{}{#1}}
\input{cyracc.def}
\font\tencyr=wncyr10
\def\cyr{\tencyr\cyracc}

\newtheorem{propo}[equation]{Proposition}
\newtheorem{corol}[equation]{Corollary}
\newtheorem{theor}[equation]{Theorem}
\newtheorem{lemma}[equation]{Lemma}
\theoremstyle{definition}
\newtheorem{defin}[equation]{Definition}

\newtheorem{assum}[equation]{Assumption}
\theoremstyle{remark}
\newtheorem{remar}[equation]{Remark}

\newtheorem{oppro}[equation]{Open Problems}

\numberwithin{equation}{section}

\newcommand{\cB }{\mathcal{B}}

\newcommand{\co }{\mathrm{co\,}}
\newcommand{\copr }{\varDelta }

\newcommand{\ndN }{\mathbb{N}}

\newcommand{\ndZ }{\mathbb{Z}}
\newcommand{\ot }{\otimes }

\newcommand{\YD }{Yetter--Drinfel'd }

\newcommand{\fil}[1]{\mathcal{F}_{#1}}
\newcommand{\filt}{\mathcal{F}}

\newcommand{\gr}{\operatorname{gr}}
\newcommand{\id}{\operatorname{id}}
\newcommand{\mbbA}{\mathbb{A}}
\newcommand{\lYDcat}[1]{{}^{#1}_{#1}\mathcal{YD}}
\newcommand{\llb}{\llbracket}
\newcommand{\rrb}{\rrbracket}
\newcommand{\db}[1]{\llb #1 \rrb}
\newcommand{\Sd}[1]{\text{\cyr X}(#1)}
\newcommand{\bn}[1]{\mathbf{#1}}
\newcommand{\mul}{\mathrm{m}}
\newcommand{\ordpr}[1]{\otimes^>_{#1}}
\newcommand{\vep}{\varepsilon}
\newcommand{\bfi}{\mathbf{k}}
\newcommand{\slwu}[1]{V^{[#1]}}
\newcommand{\slwl}[1]{V_{[#1]}}
\newcommand{\alwu}[1]{\mathcal{V}^{[\ge #1]}}
\newcommand{\alwl}[1]{\mathcal{V}_{[\ge #1]}}
\newcommand{\aalwl}[1]{\mathcal{V}^+_{[\ge #1]}}

\newcommand{\ilwu}[1]{\mathcal{I}^{[\ge #1]}}
\newcommand{\ilwl}[1]{\mathcal{I}_{[\ge #1]}}
\newcommand{\im}{\operatorname{Im}}

\newlength{\tl}

\title[On a factorization of graded Hopf algebras]{%
On a factorization of graded Hopf algebras using Lyndon words}
\author{M.~Gra\~na}
\thanks{M.\,G.~was partially supported by CONICET and Fundaci\'on Antorchas}
\thanks{M.\,G.~thanks for the warm
hospitality of Mathematisches Institut der Universit\"at Leipzig, where the work on this
paper begun.}

\address{Mat\'\i as Gra\~na, Depto de Matem\'atica - FCEyN,
Universidad de Buenos Aires, Pab I - Ciudad Universitaria
(1428) Buenos Aires - Argentina}
\email{matiasg@dm.uba.ar}

\author{I.~Heckenberger}
\address{Istv\'an Heckenberger, Mathematisches Institut, Universit\"at Leipzig,
         Augustusplatz 10-11, 04109 Leipzig, Germany.}
\email{heckenberger@math.uni-leipzig.de}

\begin{document}

\begin{abstract}
%	This is draft. Do not read it. You have been warned.
 We find a generalization of the restricted PBW basis for pointed Hopf algebras over abelian groups
constructed by Kharchenko. We obtain a factorization of the Hilbert series for a wide class of
graded Hopf algebras. These factors are parametrized by Lyndon words, and they are the Hilbert series
of certain graded Hopf algebras.


\end{abstract}
\maketitle

\section{Introduction}

Hopf algebras \cite{b-Sweedler} are far from being classified. Up to now there
are two main directions of study: semisimple and pointed Hopf algebras. This
paper is mainly a contribution to the latter, although all considerations are
performed in a more general context. Specifically, we work with Hopf algebras
$H$ generated by a Hopf subalgebra $H_0$ and a vector space $V$ satisfying
properties \eqref{eq:haf1} and \eqref{eq:haf2} below. This includes in
particular pointed Hopf algebras generated by group-like and skew-primitive
elements.  In Theorem~\ref{th:isom} we prove a factorization result about the
Hilbert series of $\gr \,H$.  Moreover, with Theorem~\ref{th:prna} we show that
to each factor one can associate in a natural way a graded Hopf algebra which
projects onto a Nichols algebra. These are the main results of the present
paper.

Kharchenko \cite{a-Khar99} proved that when $H_0$ is the group algebra of an abelian
group and it acts on $V$ by characters, $H$ admits a
restricted PBW basis. The PBW generators in this basis are labelled by Lyndon words on an alphabet
given by a set of skew-primitive elements. Examples, where $H_0$ is the group algebra of a nonabelian
group \cite{inp-MilSchn00} \cite{inp-Grana99} \cite{a-AndrGr03},
indicate that in general one can not expect Kharchenko's result to hold in its strong form.
Nevertheless, by extending his ideas we were able to construct a basis of $\gr \,H$ using ordered
products of subquotients of it.  In the particular case of Kharchenko's setting one recovers the PBW
basis. This is possible because graded braided Hopf algebras generated by one primitive element are
easy to classify. The difficulty in the more general setting arises from the fact that the structure
of graded braided Hopf algebras generated by an irreducible \YD module over $H_0$ is not known.

A generalization of Kharchenko's PBW theorem in a different direction was done by Ufer
\cite{a-Ufer03}. Instead of Hopf algebras of diagonal type he considers Hopf algebras with
triangular braidings. On the one hand Ufer is able to give a restricted PBW basis. On the other hand
some information about the relations of the Hopf algebra is lost. Although we believe that it is
possible to obtain a generalization of Ufer's approach to our context, we stick to a simpler setting
for two reasons.  First, valuable additional information can be obtained in our starting context.
Second, the proofs in the triangular case would be even more technical, obscuring the essential
arguments.

The proof of the main results of the present paper was possible due to taking advantage both from
the lexicographic and the inverse lexicographic order on the set of monotonic super-words built from
an alphabet of Lyndon words. This leads in a natural way to the construction of subquotients of a
graded Hopf algebra. In this way Kharchenko's PBW theorem becomes more transparent. Note that Ufer's
more technical proof stems from the fact that in his setting the inverse lexicographic order on the
set of monotonic super-words can not be used.

Kharchenko's PBW theorem turned out to be essential in the construction of the Weyl
groupoid \cite{a-Heck04c} corresponding to a Nichols algebra of diagonal type. This groupoid played
the crucial role in the classification of such Nichols algebras \cite{a-Heck04e}. In turn, the
knowledge of these Nichols algebras is important for example for the lifting
method of Andruskiewitsch and Schneider \cite{a-AndrSchn98} to classify pointed Hopf algebras.
We consider the factorization
theorem in this paper to be an important step towards the generalization of the Weyl groupoid to a
wider class of Nichols algebras.

In this paper $\bfi $ is an arbitrary field, and all algebras have $\bfi $ as their base field.
The symbol $\ot $ refers to tensor product over $\bfi $.
We will write $\mul $, $\copr $, $\varepsilon $ and $S$ for the product, coproduct, counit, and the antipode
of a Hopf algebra.

We conclude the introduction with two general results about Hopf algebras, for which we did not find
references in the literature.

\begin{propo}\label{pr:whab}
	Let $A=(A,\mul,\copr)$ be a bialgebra, $B\subseteq A$ a subalgebra and $I\subseteq B\cap\ker(\varepsilon)$ an
	ideal of $B$. Suppose furthermore that
	\begin{equation}\label{eq:abi}
		\copr(B)\subseteq B\ot B+I\ot A,\quad \copr(I)\subseteq B\ot I+I\ot A.
	\end{equation}
	Then the bialgebra structure on $A$ induces a bialgebra structure on $B/I$.
\end{propo}
\begin{proof}
	Notice first that $B/I$ is an algebra. Let us take $\bar{\copr}:B\to A/I\ot A/I$ as
	$\bar{\copr}=(\pi\ot\pi)\copr i$, for $\pi:A\to A/I$ the projection and $i:B\to A$ the inclusion.
	By using the first formula in \eqref{eq:abi}, one gets that
	$\bar{\copr}(B)\subseteq B/I\ot B/I$.
	By using the second formula in \eqref{eq:abi}, one gets that $\bar{\copr}(I)=0$. Then
	$\bar{\copr}$ induces a map $\tilde{\copr}:B/I\to B/I\ot B/I$. 
	Also, $\varepsilon$ induces a map $\tilde{\varepsilon}:B/I\to\bfi$. It is clear that
	$\tilde{\copr}$ is coassociative and $\tilde{\varepsilon}$ is a counit for $\tilde{\copr}$,
	whence $B/I$ is a coalgebra. It is immediate that $\tilde{\varepsilon}$ is an algebra
	map. We must prove then that $\tilde{\copr}$ is an algebra map. For $a,b\in B$ we compute
	\begin{align*}
		\tilde{\copr}(\pi(a)\pi(b)) &= \tilde{\copr}(\pi(ab)) = (\pi\ot\pi)\copr(ab)
				= (\pi\ot\pi)(a_{(1)}b_{(1)}\ot a_{(2)}b_{(2)}) \\
			&= \pi(a_{(1)}b_{(1)}) \ot \pi(a_{(2)}b_{(2)})
				= \pi(a_{(1)})\pi(b_{(1)}) \ot \pi(a_{(2)}b_{(2)}) \\
			&= \pi(a_{(1)})\pi(b_{(1)}) \ot \pi(a_{(2)})\pi(b_{(2)}).
	\end{align*}
	In the first equality above we used that $\pi|_B:B\to B/I$ is an algebra map, and the fifth one
	is obtained from $a_{(1)}b_{(1)}\ot a_{(2)}b_{(2)}\in B\ot A$. The last equality holds
	by $\pi(a_{(1)})\pi(b_{(1)})\ot a_{(2)}b_{(2)}\in B/I\ot B$.
\end{proof}

For the notion of braided Hopf algebras one may consult for example \cite{inp-Takeuchi00}.

\begin{propo}\label{pr:pNa}
Let $\pi :R\to T$ be a surjective map of braided graded Hopf algebras which is an isomorphism in
degree $0$, and assume that $R$, $T$ are finite dimensional in each degree.
Then the quotient $\eta (R,t)/\eta (T,t)$ of Hilbert series is again a series with nonnegative
integer coefficients.
\end{propo}

\begin{proof}
	Since $\ker \pi $ is homogeneous, the spaces $(\ker \pi )^i/(\ker \pi )^{i+1}$
	inherit the grading of $R$. Set
	$$
		R'=\mathop{\oplus }_{i=0}^\infty (\ker \pi )^i/(\ker \pi )^{i+1},\text{ where }
		(\ker \pi )^0/(\ker \pi )^1=R/(\ker \pi )\simeq T.
	$$
	Then $R$ and $R'$ are isomorphic as graded vector spaces. Moreover, $R'$ is a braided Hopf
	algebra, and $T\subset R$ is a braided Hopf subalgebra with a projection $\pi ':R'\to T$.
	Therefore the braided version of the Radford-Majid theorem gives that
	$R'\simeq R'{}^{\co T}\ot T$ as graded vector spaces. This implies that
	$\eta (R,t)=\eta (R',t)=\eta (R'{}^{\co T},t)\eta (T,t)$.
\end{proof}

\section{Lyndon words}
Let $A=\{\bn 1,\bn2,\dots,\bn d\}$ be a totally ordered set by $\bn 1<\bn 2<\cdots<\bn d$. We think of $A$ as
an \emph{alphabet} and $\bn 1,\cdots,\bn d$ as the \emph{letters} of $A$. Let $\mbbA$ be the set of
non-empty words in this alphabet, and let $\emptyset$ denote the empty word. For a word
$u=a_1a_2\cdots a_r$ with $a_i\in A$ $\forall i$, we say that $r$ is the \emph{length} of $u$ and we
write $r=|u|$.  We consider on $\mbbA$ the lexicographic order $<$. This means that $u<v$ if and
only if $v=uu'$ for some $u'\in\mbbA$ or $u=w\bn iu'$ and $v=w\bn jv'$, where $\bn i<\bn j$ and
$w,u',v'\in\{\emptyset\}\cup\mbbA$.

A word $u\in\mbbA$ is called a \emph{Lyndon word} if $u=u_1u_2$ with $u_1,u_2\in\mbbA$ implies that
$u<u_2$. For example: letters are Lyndon words, $\bn i\bn j$ is a Lyndon word for $\bn i<\bn j$,
$\bn 1\bn 2\bn 1\bn 2\bn 2$ is a Lyndon word, and $\bn 1\bn 2\bn 1\bn 2$ is not a Lyndon word.  We
write $L=\{u\in\mbbA\ |\ u\text{ is a Lyndon word}\}$.

\begin{propo}\cite[Prop.~5.1.3]{b-Loth83}\label{pr:sd}
	A word $u$ is a Lyndon word if and only if $u\in A$ or $u=vw$ with $v,w\in L$ and
	$v<w$. More precisely, if $w$ is the proper right factor of maximal length of $u=vw\in L$
	that belongs to $L$, then also $v\in L$ and $v<vw<w$.\hfill\qed
\end{propo}
If $u\in L$, $|u|\ge 2$, then the decomposition $u=vw$ in Prop.~\ref{pr:sd} with $w$ of maximal
length is called the \emph{Shirshov decomposition} of $u$. We write $\Sd u=(v,w)$.

\begin{lemma}\label{lm:nsd}
	Let $u,v\in L$, $u<v$, $|u|\ge 2$, and let $\Sd u=(u_1,u_2)$. Then exactly one of the following possibilities
	occurs:
	\begin{enumerate}
		\item $\Sd{uv}=(u,v)$ and $u_2\ge v$.
		\item $\Sd{uv}=(u_1',u_1''u_2v)$ for some words
			$u_1',u_1''$ such that $u_1=u_1'u_1''$ and $u_2<v$ (here $u_1''$ may be empty).
	\end{enumerate}
\end{lemma}
\begin{proof}
	This is equivalent to \cite[Prop. 5.1.4]{b-Loth83}
\end{proof}

\begin{lemma}\cite[Lemma 4]{a-Khar99}
	\label{lm:kha4}
	Let $u,v\in L$, $u=u_1u_2$ for $u_1,u_2\in\mbbA$ and suppose that $u_2<v$.
	Then $uv<u_1v$.\hfill\qed
\end{lemma}

We take on $L$ the lexicographic order. Thus, $L$ is a new alphabet containing the original alphabet
$A$, and following Kharchenko \cite{a-Khar99}
we say that the elements of $L$ are \emph{super-letters}. Words in super-letters are
called \emph{super-words}. The \emph{length} $|w|$ of a super-word $w$ is the sum of the lengths of its
super-letters. A \emph{monotonic super-word} is a non-increasing word on the set of
super-letters, i.e., a (possibly empty) word $v_1\dots v_n$ such that $v_i\in L$ and
$v_1\ge v_2\ge \dots\ge v_n$. Let $M$ denote the set of monotonic super-words.
In what follows the notation $v_1\cdots v_n\in M$ will mean that $v_1,\ldots ,v_n\in L$ and
$v_1\ge \cdots \ge v_n$. Sometimes we also write $v_1^{m_1}\cdots v_n^{m_n}\in M$, in which case
we mean that $v_1,\ldots ,v_n\in L$, $v_1>\cdots >v_n$, and $m_1,\ldots ,m_n\ge 1$.
Monotonic super-words are lexicographically ordered on the alphabet
of super-letters.  Notice that the empty super-word is the smallest super-word.
For a super-letter $u$, we shall write
\begin{equation}\label{eq:mbu}
	\begin{split}
		&L_{>u}=\{v\in L\;|\; v>u\},\quad L_{\ge u}=\{v\in L\;|\; v\ge u\}, \\
		&M_{>u}=\{v_1\cdots v_r\in M\;|\; r\ge 1,\ v_1\ge \dots \ge v_r>u\}, \\
		&M_{\ge u}=\{v_1\cdots v_r\in M\;|\; r\ge 1,\ v_1\ge \dots \ge v_r\ge u\}.
	\end{split}
\end{equation}

\begin{theor}(Lyndon, see \cite[Thm. 5.1.5 and Prop. 5.1.6]{b-Loth83})\label{th:lyn}
	A word in $\mbbA$ can be written in a unique way as a monotonic super-word.
	Moreover, if $u=v_1\dots v_n\in M$ then $v_n$ is the smallest
	right factor of $u$ (smallest with respect to lexicographic order in $\mbbA $).\hfill\qed
\end{theor}

As an example, the word $\bn 1\bn 2\bn 3\bn 1\bn 2\bn 3\bn 3\bn 1\bn 2\bn 2\bn 1\bn 2\bn 3$ is decomposed as a
monotonic super-word as
$$\bn 1\bn 2\bn 3\bn 1\bn 2\bn 3\bn 3\bn 1\bn 2\bn 2\bn 1\bn 2\bn 3
	=(\bn 1\bn 2\bn 3\bn 1\bn 2\bn 3\bn 3)(\bn 1\bn 2\bn 2\bn 1\bn 2\bn 3),$$
and in turn, 
$\Sd{\bn 1\bn 2\bn 3\bn 1\bn 2\bn 3\bn 3}=(\bn 1\bn 2\bn 3,\bn 1\bn 2\bn 3\bn 3)$. %=((\bn 1)(\bn 2\bn 3))((\bn 1)((\bn 2\bn 3)(\bn 3))).$$

\begin{lemma}\label{lm:dbig}
	Let $w$ be a super-word. Then the decomposition of $w$ as a monotonic super-word is $\ge w$
	with respect to the lexicographic order on super-words.
\end{lemma}
\begin{proof}
We proceed by induction on the number of super-letters of $w$.
If $w$ is a super-letter then we are done. Otherwise $w=w_1\cdots w_n$, where $n\ge 2$ and $w_i\in L$ $\forall i$.
Again, if $w_1\ge w_2\ge \ldots \ge w_n$ then we are done. On the other hand, if $w_i<w_{i+1}$
for some $i$, then $w'_i:=w_iw_{i+1}\in L$ by Proposition~\ref{pr:sd},
and one has $w=w_1\cdots w_{i-1}w'_iw_{i+2}\cdots w_n$ with $w'_i>w_i$.
Hence the claim follows from the induction hypothesis.
\end{proof}

\begin{lemma}\cite[Lemma 5]{a-Khar99}\label{lm:khal5}
	Let $w=w_1\cdots w_m$ and $v=v_1\cdots v_n$ be monotonic super-words. Then $w<v$ (with respect to the
	lexicographic order on $M$) if and only if $w_1\cdots w_m<v_1\cdots v_n$
	with respect to the lexicographic order on $\mbbA $. \hfill \qed
\end{lemma}

The following technical lemma will be needed in the proof of Theorem~\ref{th:isom}.

\begin{lemma}\label{lm:wdec}
	Let $w=w_1\cdots w_m$ and $v=v_1\cdots v_n$ be nonempty monotonic super-words with $w\ge v$
	and assume that $v_1=\cdots =v_r>v_{r+1}$, where $r\le n$. For all $i\le m$ let $(w'_i,w''_i)\in M\times M$
	such that either $w_i'\ge w_i$ or $w'_i=\emptyset $, $w''_i=w_i$. Then
        the pair $(w'_1\cdots w'_m,w''_1\cdots w''_m)$
	satisfies one of the following relations.
	\begin{enumerate}
		\item \label{lm:wdec1}
			$w'_1\cdots w'_m>v_{r+1}\cdots v_n$,
		\item \label{lm:wdec2}
			$w''_1\cdots w''_m>v_1^r$,
		\item \label{lm:wdec3}
			$(w'_1\cdots w'_m,w''_1\cdots w''_m)=(v_{r+1}\cdots v_n,v_1^r)$.
	\end{enumerate}
	Moreover, if $w>v$ then \eqref{lm:wdec3} can not hold.
\end{lemma}

\begin{proof}
	Assume first that $w_1>v_1$. Then either $w'_1>v_1>v_{r+1}\cdots v_n$ or $w'_1=\emptyset $, $w''_1=w_1>v_1$.
	In the first case we have relation \eqref{lm:wdec1} and in the second one relation \eqref{lm:wdec2} is fulfilled.
	On the other hand, if $w_1\le v_1$ then because of $w\ge v$ and $w$ is monotonic, one has $m\ge r$ and $w_i=v_1$
	for all $i\le r$. Suppose first that there exists $i\le r$ such that $w'_i\not=\emptyset $ and $w'_j=\emptyset $
	for all $j<i$. In this case $w'_i\ge v_1>v_{r+1}\cdots v_n$ and hence \eqref{lm:wdec1} holds.
	It remains to consider the case $w'_i=\emptyset $ for all $i\le r$. Then one has $w''_1\cdots w''_r=v_1^r$.
	Therefore, if $w''_i\not= \emptyset $ for some $i>r$ then again relation \eqref{lm:wdec2} holds.
	Otherwise $w'_i\ge w_i$ for all $i>r$. Then one has either
	$w'_1\cdots w'_m=w'_{r+1}\cdots w'_m>v_{r+1}\cdots v_n$, in which case \eqref{lm:wdec1} holds,
        or $m=n$, $w=v$, and $w'_i=w_i$ for all $i>r$.
	The latter relations imply \eqref{lm:wdec3}.
\end{proof}

Let $H_0$ be a Hopf algebra and let $V=\oplus_{\bn i=\bn 1}^{\bn d}V_{\bn i}$ be a direct sum of
\YD modules over
$H_0$.  Let $TV$ be the tensor algebra of $V$. For simplicity, we will omit the $\ot$ symbol in
roducts of elements of $TV$.
Let $.$ and $\delta $ denote the left action and the left coaction of $H_0$ on $TV$, respectively.
We will use Sweedler notation $\delta (x)=x_{(-1)}\ot x_{(0)}$ for $x\in TV$.
The algebra $TV$ has a braiding $c:TV\ot TV\to TV\ot TV$.
Note that one has
\begin{align*}
 c(x\ot y)=&x_{(-1)}.y\ot x_{(0)},&
 c^{-1}(x\ot y)=&y_{(0)}\ot S^{-1}(y_{(-1)}).x
\end{align*}
for all $x,y\in TV$.
We define
$$
	[x,y]=xy-\mul c^{-1}(x\ot y),\qquad \llb x,y \rrb = xy-\mul c(x\ot y).
$$
where $\mul$ is the multiplication in $TV$.

\begin{defin}\label{df:uvec}
	Let $a_1,\dots,a_m\in A$ and $u=a_1\dots a_m\in\mbbA$.  We write
	$V^u=V_{a_1}V_{a_2}\cdots V_{a_m}\simeq V_{a_1}\ot V_{a_2}\ot\dots\ot V_{a_m}$. The elements in
	$V^u$ will be called \emph{$u$-vectors}. If $x_u\in V^u$ and $u=vw$, then we write
	$x_u=x_v\ot x_w\in V^v\ot V^w$ (which is in general a sum of tensors)
        using the canonical isomorphism $V^u\simeq V^v\otimes V^w$.
\end{defin}

We shall inductively define bracket operations $[\;]:\oplus_{n\ge 0}V^{\ot n}\to TV$ and
$\llb\;\rrb:\oplus_{n\ge 0}V^{\ot n}\to TV$ as follows. Let $x_u$ be a $u$-vector.
\begin{enumerate}
	\item If $u$ has length $0$ or $1$, then $[x_u]=\db{x_u}=x_u$.
	\item If the word $u$ is a Lyndon word and $\Sd u=(v,w)$, then
		$[x_u] = [[x_v]\,,\,[x_w] ]$ and $\db{x_u} = \db{\db{x_v}\,,\,\db{x_w}}$
		(see Def.~\ref{df:uvec}).
	\item If the word $u$ is decomposed as a monotonic super-word by $u=v_1\cdots v_r$, then
		$[x] = [x_{v_1}]\cdot[x_{v_2}]\cdots[x_{v_r}]$ and
		$\db x = \db{x_{v_1}}\cdot\db{x_{v_2}}\cdots\db{x_{v_r}}$.
\end{enumerate}

\begin{remar}
Recall that the braided antipode $S_{TV}$ of the braided Hopf algebra $TV$ satisfies
$S_{TV}(x)=x_{(-1)}S(x_{(0)})$, where $S$ is the antipode of $TV\#H_0$.
Moreover, $S_{TV}(xy)=\mul c(S_{TV}(x)\ot S_{TV}(y)) =(x_{(-1)}.S_{TV}(y))S_{TV}(x_{(0)})$. With
these formulas it is easy to see that for any $u\in L$ and any $u$-vector $x$ one has
$S_{TV}([x])=(-1)^{|u|-1}\db{x}$. Therefore most of the following considerations can be performed
without difficulties with $\db{\,}$'s instead of $[\,]$'s. Even if $\db{\,}$'s seem to be more
natural, we will follow the tradition of Kharchenko \cite{a-Khar99} with $[\,]$'s.
\end{remar}

\begin{defin}
	If $u$ is a Lyndon word and $x$ is a $u$-vector, then $[x]$ will be called \emph{a
	$u$-$[\,]$-letter}.
	If $u$ is a monotonic super-word and $x$ is a $u$-vector then $[x]$ will be called
	\emph{a $u$-$[\,]$-word}.
	Let $V^{[u]}=[V^u]$ denote the space of $u$-$[\,]$-words.
	A $[\,]$-letter ($[\,]$-word) is a $u$-$[\,]$-letter ($u$-$[\,]$-word) for some super-letter
	(monotonic super-word) $u$.
\end{defin}

\begin{lemma}\label{lm:subskew}
%       Let $u$ be a super-word.
        Let $u$ be a monotonic super-word.
        For $x,y\in TV$, let $c^{-1}(x\ot y)=:\sum {}^xy\ot x^y$. Then
	\begin{enumerate}
		\item\label{lm:subskew1}
			if $x$ is a $u$-vector and $h\in H_0$ then one has $[h.x]=h.[x]$ and
			$[x]_{(-1)}\ot [x]_{(0)}=x_{(-1)}\ot [x_{(0)}]$;
		\item\label{lm:subskew1b}
			if $x$ is a $u$-vector and $y$ is a $v$-vector for a monotonic super-word $v$, then
			${}^{[x]}y \ot [x]^y = {}^xy \ot [x^y]$ and ${}^x[y]\ot x^{[y]} = [{}^xy] \ot x^y$;
		\item\label{lm:subskew2} for all $x,y,z\in TV$,
			\begin{align*}
				[x,yz] &=[x,y]\; z + \sum {}^xy\;[x^y,z];
			\end{align*}
		\item\label{lm:subskew3} for all $x,y,z\in TV$,
			\begin{align*}
				[ [x,y],z ] - [ x, [y,z] ] &= \sum [x,{}^yz] y^z - \sum {}^xy[x^y,z].
			\end{align*}
	\end{enumerate}
\end{lemma}
\begin{proof}
	\eqref{lm:subskew1} follows from the fact that the brading is a map of \YD modules.
	\eqref{lm:subskew1b} follows from \eqref{lm:subskew1} and the braid relation.
	\eqref{lm:subskew2} and \eqref{lm:subskew3} are straightforward calculations using the definition of
	$[\,]$ and the braid relation for $c^{-1}$.
\end{proof}

\stepcounter{equation}
We prove now a variant of \cite[Lemma 6]{a-Khar99}.
\begin{lemma}\label{lm:khal6}
	Let $X,Y$ be a $u$-$[\,]$-letter and a $v$-$[\,]$-letter respectively, where $u,v\in L$ and $u<v$.
	Then $[X,Y]$ is a homogeneous linear combination of products of $[\,]$-letters corresponding to
        super-letters in $L_{\ge uv}$.
\end{lemma}
\begin{proof}
	Let $z=uv$. If $\Sd z=(u,v)$, then $[X,Y]$ is a $z$-$[\,]$-letter and we are done.
        We proceed by induction on $|uv|$, since we already know the lemma for the case $|u|=|v|=1$.

	Let $m=|uv|$ and suppose that the lemma holds for $u',v'\in L$, $u'<v'$, $X\in V^{[u']}$,
        $Y\in V^{[v']}$, where either $|u'v'|<m$ or $|u'v'|=m$, $u'<u$
        (notice that there are only finitely many Lyndon words of a given length).
	As noted above, we may suppose that $\Sd {uv}\neq (u,v)$. Then, if $\Sd u=(u_1,u_2)$, we must have
	$u_2<v$ because of Lemma~\ref{lm:nsd}.  Let $X=[x]$ for $x\in V^{u}$, $Y=[y]$ for $y\in V^v$,
	and let $x=x_{u}=x_{u_1}\ot x_{u_2}$, $X_{u_1}\ot X_{u_2}=[x_{u_1}]\ot [x_{u_2}]$.  Thanks to
	Lemma~\ref{lm:subskew}\eqref{lm:subskew3},\eqref{lm:subskew1b}, we have
	\begin{equation}\label{eq:1}
		\begin{split}
		[X,Y]=&[[X_{u_1},X_{u_2}],Y]\\
			= &[X_{u_1},[X_{u_2},Y]] \\
				&\hspace{.5cm}+ \sum [X_{u_1},{}^{X_{u_2}}Y] (X_{u_2})^Y
				- \sum {}^{X_{u_1}}X_{u_2}[X_{u_1}^{X_{u_2}},Y] \\
			= &[X_{u_1},[X_{u_2},Y]] \\
				&\hspace{.5cm}+ \sum [ X_{u_1},[{}^{x_{u_2}}y]]\cdot[x_{u_2}^y]
				- \sum [{}^{x_{u_1}}x_{u_2}]\cdot[ [x_{u_1}^{x_{u_2}}],Y].
		\end{split}
	\end{equation}
	We start by considering the first summand on the right hand side.  By the induction hypothesis,
	$[X_{u_2},Y]$ is a homogeneous linear combination of products of $[\,]$-letters corresponding to
	super-letters in $L_{\ge u_2v}$, and the degree of these products is $|u_2v|$.  By
	Lemma~\ref{lm:subskew}\eqref{lm:subskew2}, and since $u_1<u_2<u_2v$, $[X_{u_1},[X_{u_2},Y]]$ is
	a combination of products of $[\,]$-letters corresponding to super-letters in $L_{\ge
	u_2v}\subset L_{>uv}$ on one hand, and a bracket $[h.X_{u_1},X']$ on the other hand,
	where $h\in H_0$ and $X'$ is a $w$-$[\,]$-letter with $w\ge u_2v$. By induction
	hypothesis again, the latter is a linear combination of products of $[\,]$-letters corresponding
	to super-letters in $L_{\ge u_1u_2v}=L_{\ge uv}$.

	We continue with the second and third summands in \eqref{eq:1}. Concerning the
	$u_2$-$[\,]$-letters appearing there, notice that $u_2>uv$. Indeed, $u_2>u_1u_2$ and, since
	$u_1u_2$ is not the beginning of $u_2$, then we still have $u_2>u_1u_2v=uv$.
	On the other factors, which are brackets between $u_1$-$[\,]$-letters and $v$-$[\,]$-letters,
	we can apply the induction hypothesis since $|u_1v|<|uv|$ and $u_1<u_1u_2<u_2<v$. These factors
	are then linear combinations of products of $[\,]$-letters corresponding to super-letters in $L_{\ge u_1v}$,
	and $u_1v>uv$ by Lemma~\ref{lm:kha4}.
\end{proof}

\begin{lemma}
	\label{lm:etim}
	Let $u\in L$.  Any product of $[\,]$-letters corresponding to super-letters in $L_{\ge u}$ is a
	linear combination of (monotonic) $[\,]$-words corresponding to super-words in $M_{\ge u}$.
\end{lemma}
\begin{proof}
	Let $t=v_1\cdot v_2\cdots v_n$ be a super-word, where $v_i\in L_{\ge u}$ $\forall i$. Let
	$x=x_t=x_{v_1}\cdots x_{v_n}$ be a $t$-vector, where here $t$ is considered as a word on $A$. We
	call $x_{[t]}:=[x_{v_1}]\cdots[x_{v_n}]$ \emph{a $t$-$[\,]$-vector} (notice that the term
	``$t$-$[\,]$-word'' is reserved for when $t$ is monotonic).
	We proceed by induction. We take on the set of super-words the lexicographic order $\prec$ given
	by the order $<$ on $L$.  Suppose that the result is true for $w$-$[\,]$-vectors when $|w|<|t|$
	or when $|w|=|t|$ and $w\succ t$.  If $v_1\ge v_2\ge \dots\ge v_n$, there is nothing to prove.
	Otherwise, let $i$ be such that $1\le i<n$ and $v_i<v_{i+1}$. Let
	$t'=v_1\dots v_{i+1}v_i \dots v_n$ and $t''=v_1\dots (v_iv_{i+1})\dots v_n$, where the factor
	$(v_iv_{i+1})$ stands for the Lyndon word $v_iv_{i+1}$. Let
	\begin{align*}
		x_{[t']} &:=[x_{v_1}]\ot\dots\ot c^{-1}([x_{v_i}]\ot[x_{v_{i+1}}])\ot\dots\ot[x_{v_n}], \\
		x_{[t'']}&:=[x_{v_1}]\ot\dots\ot [[x_{v_i}], [x_{v_{i+1}}]]\ot\dots\ot[x_{v_n}].
	\end{align*}
	Notice that $x_{[t']}$ is a $t'$-$[\,]$-vector where $t'$ is a super-word in the same
	super-letters as $t$, and $t'\succ t$. Also, by Lemma~\ref{lm:khal6}, and since
	$v_iv_{i+1}>v_i$, $x_{[t'']}$ is a linear combination of $w$-$[\,]$-vectors, where $w$ runs on
	super-words $\succ t$ having only super-letters in $L_{\ge u}$. The last thing to notice is that
	$x_{[t]}=x_{[t']}+x_{[t'']}$. Therefore, the induction hypothesis implies the claim.
\end{proof}

\begin{defin}
	Let $x\in TV\setminus\{0\}$, $x=\sum_nx_n$, where $x_n\in V^{\ot n}$. The greatest $n$ such that
	$x_n\neq 0$ will be called the \emph{degree} of $x$. Let $m$ be the degree of $x$ and let
	$x_m=\sum_{u\in\mbbA_m}x_u$, where $\mbbA_m$ is the set of words of length $m$ and $x_u$ is a
	$u$-vector for all $u$. Let $v$ be the least $u$ (with respect to the order in $\mbbA_m$) such
	that $x_u\neq 0$. Then $x_v$ will be called \emph{the leading vector} of $x$.
\end{defin}

\begin{lemma}
	\label{lm:coeff1}
	Let $x\in TV$ be a nonzero $u$-vector for a monotonic super-word $u$. Then $x$ is
	the leading vector of $[x]$.
\end{lemma}
\begin{proof}
	We prove first the lemma for Lyndon words $u$ by induction on
	$|u|$.  If $|u|$ is $0$ or $1$, the result is clear. If $\Sd u=(v,w)$, then $[x]=\left[ [x_v],[x_w]
	\right]=(\mul-\mul c^{-1})([x_v]\ot[x_w])$. By the induction hypothesis, $[x_v][x_w]$ is a sum
	of products of $v'$-vectors and $w'$-vectors, where $v'$ runs on words $\ge v$ and $w'$ runs on
	words $\ge w$. For such $v',w'$ we have $u=vw\le v'w'$, whence $[x_v][x_w]$ is a sum of
	$u'$-vectors with $u'\ge u$.  Furthermore, the equality holds if and only if $v'=v$ and $w'=w$.
	Thus, by using induction hypothesis, the leading vector of $[x_v][x_w]$ is $x_vx_w$.  For the
	term $\mul c^{-1}([x_v]\ot[x_w])$ we reason similarly: as $u$ is a Lyndon word, $u=vw<w<wv$,
	from where $\mul c^{-1}([x_v]\ot[x_w])$ is a sum of $u'$-vectors where $u'$ runs on words $\ge
	wv>u$. Therefore, such $u'$-vectors do not contribute to the leading vector of $[x_vx_w]$.

	If $u$ is not a super-letter, then by Theorem~\ref{th:lyn} we have $u=v_1\cdots v_n\in M$.
	By definition $[x]=[x_{v_1}]\cdots[x_{v_n}]$, which, by the
	previous step, is a sum of $v'_1\cdots v'_n$-vectors with $v'_i\ge v_i$ $\forall i$, and its
	leading vector is $x_{v_1}\cdots x_{v_n}=x$.
\end{proof}

Recall that $M$ is the set of monotonic super-words.
\begin{corol}\label{co:dsd}
	One has $TV=\bigoplus_{u\in M}V^{[u]}$.
\end{corol}
\begin{proof}
	Since letters in $A$ are also super-letters, the spaces $V^{[u]}$ generate $TV$ thanks to
        Lemma~\ref{lm:etim}. The linear independence of these spaces follows immediately from
        Lemma~\ref{lm:coeff1}.
\end{proof}

\begin{corol}\label{co:debr}
	Let $u\in \mbbA $ and $x$ be a $u$-vector. Then $x-[x]$
	is a linear combination of $w$-$[\,]$-words with $w>u$.
\end{corol}

\begin{proof}
	This assertion can be proven with the help of bases of $\oplus _{|v|=|u|}V^v$ and $\oplus _{|v|=|u|}V^{[v]}$
	which are obtained from each other using triangular matrices.
	Alternatively, Lemmas~\ref{lm:khal5} and \ref{lm:coeff1} imply
	that $x-[x]$ is a linear combination of $w$-vectors
	with $w>u$, and proceed by induction.
\end{proof}

Corollary~\ref{co:debr} allows us to give a description of products of $[\,]$-letters, which is
different from the one in Lemma~\ref{lm:etim}.

\begin{corol}\label{co:etim2}
	Let $n\in \ndN $, $v_1,\ldots ,v_n\in L$ and let $X_i$ be a $v_i$-$[\,]$-letter
	$\forall i$.  Then $X_1\cdots X_n$ is a linear combination of
	$w$-$[\,]$-words, where $w$ runs over monotonic super-words
	$\ge v'_1\cdots v'_{n'}$, and $v'_1\cdots v'_{n'}$ is the decomposition of $v_1\cdots v_n$
	as a monotonic super-word. \hfill \qed
\end{corol}

\begin{defin}\label{df:op}
	Let $K$ be a totally ordered set and let there be vector spaces $W_k$ for each $k\in K$.
	We define $\ordpr{k\in K}W_k$ to be the direct sum of vector spaces
	$W_{k_1}\ot \dots \ot W_{k_r}$ where $k_1>k_2>\dots>k_r$.
\end{defin}

For a vector space $W$, we write $T^+W=\oplus_{n\ge 1}W^{\ot n}$.  Notice that for each $u\in L$,
$T^+V^{[u]}$ is a non-unital subalgebra of $T^+V$ in the category $\lYDcat{H_0}$.
Further, the $\ndZ $-grading of $TV$ induces a $\ndZ $-grading on its subalgebras $TV^{[u]}$.
In particular, the Hilbert series $\eta (T\slwu u,t)$ of $T\slwu u$
is a series in the variable $t^{|u|}$.

\begin{theor}\label{th:md}
	One has
	$$
	T^+V\simeq\ordpr {u\in L}T^+V^{[u]}.
	$$
	More precisely, the map $\mu:\ordpr {u\in L}T^+V^{[u]}\to T^+V$, which is the
	multiplication map in each summand, is an isomorphism in the category $\lYDcat{H_0}$.
	In particular, the Hilbert series of $TV$ is
	$\eta(TV,t)=\prod_{u\in L}\eta(T\slwu u,t)$.
\end{theor}
\begin{proof}
	This is a reformulation of Corollary~\ref{co:dsd}.
\end{proof}

\section{Hopf algebras generated by a Hopf subalgebra and a vector space}
Our aim is now to apply the results of the previous section to arbitrary Hopf algebras.
In many (though not all) cases we will be able to provide some new structure results.

Let $H$ be a Hopf algebra with a filtration
$$
0=\fil{-1} \subset \fil 0 \subset \fil 1 \subset \fil 2 \subset \dots \subseteq H.
$$
The filtration is a \emph{Hopf algebra filtration} if
\begin{enumerate}
	\item It is a filtration: $H=\bigcup_{n\in\ndN}\fil n$,
	\item it is an algebra filtration: $\fil n\fil m\subseteq\fil{n+m}$,\label{af}
	\item it is a coalgebra filtration: $\copr(\fil n)\subseteq\sum_{i+j=n}\fil i \otimes \fil j$, and
	\item it behaves well with the antipode: $S(\fil n)\subseteq\fil n$.
\end{enumerate}
Notice that in this case $\fil 0$ is a Hopf subalgebra of $H$.
We will consider Hopf algebras with a Hopf algebra filtration satisfying a stronger
version of condition \eqref{af}:
\begin{assum}\label{as:fgr}
	$\fil n\fil m=\fil{n+m}$.
\end{assum}

A Hopf algebra with a Hopf algebra filtration satisfying Assumption~\ref{as:fgr} can be presented in
the following way.

Suppose that the Hopf algebra $H$ is generated (as an algebra) by a Hopf subalgebra
$H_0$ and a vector space $V$, such that
\begin{align}
\label{eq:haf1}
	&\copr V\subseteq V\ot H_0+H_0\ot V+H_0\ot H_0,\text{ and }\\
\label{eq:haf2}
	&S(V) \subseteq KVK.
\end{align}
We define $\fil 0=H_0$, $\fil 1=H_0+H_0VH_0$, and $\fil n=(\fil 1)^n$. Then $\fil *$ is a Hopf algebra
filtration which satisfies Assumption~\ref{as:fgr}.

As an important example, assume that  $H$ is generated by grouplikes and skew-primitive
elements and take $H_0$ to be the subalgebra generated by the grouplikes and $V$ to be the subspace
generated by skew-primitives.

%\begin{question}
%	Is it true that if $H$ is pointed, $\fil 0$ is the coradical and $\fil 1$ is the space
%	generated by skew primitives and $1$, then $\fil *$ coincides with the coradical
%	filtration?
%\end{question}

For $n\ge 0$ let $H_n=\fil n/\fil {n-1}$. Let then
$$
\gr_{\filt}H=\mathop{\oplus}_{n\ge 0}\fil n/\fil {n-1}
	=\mathop{\oplus}_{n\ge 0}H_n
$$
be the associated graded Hopf algebra.  Then we can consider the projection
$\pi:\gr_{\filt}H\to H_0$ and the inclusion $\iota:H_0\to\gr_{\filt}H$, and this allows to
write $\gr_{\filt}H=R\# H_0$ as the smash product of $H_0$ and the right coinvariants
$R=\{x\in \gr_{\filt}H\;|\;(\id\ot\pi)\copr(x)=x\ot 1\}$.
This is now a standard procedure: the last part is due to Radford \cite{a-Radford85} and
Majid \cite{a-Majid94a}, while the first part is a modification of the one due to
Andruskiewitsch and Schneider \cite{a-AndrSchn98}.

Also, $R$ is a braided Hopf algebra in $\lYDcat{H_0}$, and $R=\oplus_{n\ge 0}R_n$, $R_0$ is the
base field and $R_1$ is a \YD module which generates $R$ because of Assumption~\ref{as:fgr}.
We let $V:=R_1$ and we make the following assumption:
\begin{assum}
	$V=\oplus_{i=1}^dV_i$ is a direct sum of \YD modules over $H_0$.
\end{assum}

\begin{remar}
If $V$ is irreducible, the methods in this paper do not yield any information on $R$.
Otherwise, $V$ has a nontrivial flag $V=V_1\supset V_2\supset \cdots \supset V_d$ of
\YD submodules over $H_0$. Consider on the vector space $R_n$ the $\ndZ $-filtration
$R_n=R_{n,n}\supset R_{n,n+1}\supset \cdots $, where $R_{n,m}=\sum _{i_1+\cdots +i_n\ge m}V_{i_1}
V_{i_2}\cdots V_{i_n}$. Then
 $$
  R':=\mathop{\oplus} _{n=0}^\infty R'_n,\quad \text{where }
  R'_n=\mathop{\oplus} _{m\ge n} R_{n,m}/R_{n,m+1},
 $$
is a graded braided Hopf algebra in the category $\lYDcat{H_0}$, and all
of the following considerations may be applied to $R'$ instead of $R$.
\end{remar}

Consider the
projection $TV\to R$. We want to study the images under this projection of the components
$T^+V^{[u]}$ appearing in Theorem~\ref{th:md}.
We begin by considering the comultiplication in $TV\# H_0$.

We write the left $H_0$-coaction on $r\in R$ by $\delta(r)=r_{(-1)}\ot r_{(0)}$.
Recall that the coproduct in the smash product $R\# H_0$ is given by
$\copr(r\# h)=(r^{(1)}\# {r^{(2)}}_{(-1)}h_{(1)}) \ot ({r^{(2)}}_{(0)}\# h_{(2)})$,
where $\copr _R(r)=:r^{(1)}\ot r^{(2)}$ is the coproduct of the braided Hopf algebra $R$.
This notation applies in particular for $R=TV$.

\begin{propo}\label{pr:com}
	Let $X\in TV\subseteq TV\# H_0$ be a $u$-$[\,]$-letter. Then the coproduct $\copr$ of
	$TV\# H_0$ satisfies
	$$
	\copr(X)=X\ot 1 + X_{(-1)}\ot X_{(0)} + \sum_i(X'_ih_i)\ot X''_i,
	$$
	where $X'_i,X''_i\in T^+V$, $h_i\in H_0$, and each $X'_i$ is a $w_i$-$[\,]$-word with $w_i\in M_{>u}$.
\end{propo}
\begin{proof}
	We proceed by induction on $|u|$.
	If $u\in A$, we get $X\in V$ and then $\copr(X)=X\ot 1+X_{(-1)}\ot X_{(0)}$, whence we are done.
	Assume now that $|u|\ge 2$, $\Sd u=(v,w)$ and $x$ is a $u$-vector. Then $x=x_vx_w$, and we write
	$X=[x]$, $Y=[x_v]$, $Z=[x_w]$. By standard computations,
	\begin{equation}\label{eq:fdif}
		\begin{split}
			\copr(X) &=\copr([ Y, Z ]) = \copr (YZ-Z_{(0)}(S^{-1}(Z_{(-1)}).Y)) \\
			&=
			(Y^{(1)}{Y^{(2)}}_{(-1)} \ot {Y^{(2)}}_{(0)})\cdot (Z^{(1)}{Z^{(2)}}_{(-1)} \ot {Z^{(2)}}_{(0)}) \\
			&\hspace{.8cm} -({Z_{(0)}}^{(1)}{{Z_{(0)}}^{(2)}}_{(-1)} \ot {{Z_{(0)}}^{(2)}}_{(0)}) \\
			&\hspace{1.6cm} \cdot \big((S^{-1}(Z_{(-1)}).Y^{(1)})S^{-1}(Z_{(-2)}){Y^{(2)}}_{(-1)}Z_{(-4)}\\
                        &\hspace{6cm}\ot {S^{-1}(Z_{(-3)}).Y^{(2)}}_{(0)}\big).
		\end{split}
	\end{equation}
	Note that $Y$ is a $v$-$[\,]$-letter and $Z$ is a $w$-$[\,]$-letter. According to the induction
	hypothesis and Lemma~\ref{lm:subskew}\eqref{lm:subskew1}, for any $h\in H_0$, $h.Y^{(1)}$ can be
	taken to be either $\vep(h)1$, $h.Y$ or a $v'$-$[\,]$-word with $v'\in M_{>v}$. Similarly,
	$h.Z^{(1)}$ can be taken to be either $\vep(h)1$, $h.Z$ or a $w'$-$[\,]$-word with
	$w'\in M_{>w}$, and ${Z_{(0)}}^{(1)}$ to be either $1$, $Z_{(0)}$ or a $w'$-$[\,]$-word with
	$w'\in M_{>w}$.

	We begin by considering the summand
	\begin{align*}
	 (Y^{(1)}{Y^{(2)}}_{(-1)} &\ot {Y^{(2)}}_{(0)})\cdot (Z^{(1)}{Z^{(2)}}_{(-1)} \ot {Z^{(2)}}_{(0)}) \\
	 &= Y^{(1)}{Y^{(2)}}_{(-1)}Z^{(1)}{Z^{(2)}}_{(-1)} \ot {Y^{(2)}}_{(0)}{Z^{(2)}}_{(0)} \\
	 &= Y^{(1)}({Y^{(2)}}_{(-2)}.Z^{(1)}){Y^{(2)}}_{(-1)}{Z^{(2)}}_{(-1)} \ot {Y^{(2)}}_{(0)}{Z^{(2)}}_{(0)}.
	\end{align*}
	We consider the summands in which $Y^{(1)}$ is a $v'$-$[\,]$-word with $v'\in M_{>v}$. Notice
	that since these $v'$ are shorter than $v$, they belong to $M_{>vw}$. Therefore, since $Z^{(1)}$ is either
        $1$ or a $t$-$[\,]$-word with $t\in M_{>vw}$, by Lemma~\ref{lm:etim} these summands satisfy the
	claim of the Proposition. The summands in which $Y^{(1)}=1$ and $Z^{(1)}\neq 1$ also satisfy the
	claim, because $w>u$. We are thus left with
	\begin{equation}\label{eq:fsm}
		Y_{(-1)}Z_{(-1)}\ot Y_{(0)}Z_{(0)} + YZ^{(1)}{Z^{(2)}}_{(-1)} \ot {Z^{(2)}}_{(0)}.
	\end{equation}

	We consider now the other summand of $\copr(X)$. By similar reasons, we are left with
	\begin{equation}\label{eq:ssm}
		\begin{split}
		&- Y_{(-1)}Z_{(-2)}\ot Z_{(0)}\big(S^{-1}(Z_{(-1)}).Y_{(0)}\big) \\
		&\hspace{.6cm} - {Z_{(0)}}^{(1)}{{Z_{(0)}}^{(2)}}_{(-1)}\big(S^{-1}(Z_{(-1)}).Y\big)S^{-1}(Z_{(-2)})Z_{(-4)}\\
                &\hspace{6cm} \ot {{Z_{(0)}}^{(2)}}_{(0)}\big(S^{-1}(Z_{(-3)}).1\big) \\
		&=- Y_{(-1)}Z_{(-2)}\ot Z_{(0)}\big(S^{-1}(Z_{(-1)}).Y_{(0)}\big) \\
		&\hspace{.6cm} - {Z_{(0)}}^{(1)}{{Z_{(0)}}^{(2)}}_{(-1)}\big(S^{-1}(Z_{(-1)}).Y\big)
                 \ot {{Z_{(0)}}^{(2)}}_{(0)} \\
		&=- Y_{(-1)}Z_{(-2)}\ot Z_{(0)}\big(S^{-1}(Z_{(-1)}).Y_{(0)}\big) \\
		&\hspace{.6cm} - {Z^{(1)}}_{(0)} \big(\big({Z^{(2)}}_{(-2)}
                 S^{-1}({Z^{(2)}}_{(-3)})S^{-1}({Z^{(1)}}_{(-1)})\big).Y\big)
				{Z^{(2)}}_{(-1)} \ot {Z^{(2)}}_{(0)} \\
		&=- Y_{(-1)}Z_{(-1)}\ot \mul c^{-1}(Y_{(0)}\ot Z_{(0)})\\
		&\hspace{5cm} - \mul c^{-1}(Y \ot {Z^{(1)}}) {Z^{(2)}}_{(-1)} \ot {Z^{(2)}}_{(0)}.
		\end{split}
	\end{equation}
	Adding \eqref{eq:fsm} and \eqref{eq:ssm} we get
	$$
		Y_{(-1)}Z_{(-1)}\ot [Y_{(0)}, Z_{(0)}] + [Y , {Z^{(1)}}] {Z^{(2)}}_{(-1)} \ot {Z^{(2)}}_{(0)}.
	$$
	Notice that $Y_{(-1)}Z_{(-1)}\ot [Y_{(0)}, Z_{(0)}]=X_{(-1)}\ot X_{(0)}$.
	Also, when in the second summand we put $Z^{(1)}=Z$, we get $X\ot 1$. Finally, when $Z^{(1)}$ is a
	$w'$-$[\,]$-word with $w'\in M_{>w}$, by using Lemmas~\ref{lm:subskew}\eqref{lm:subskew2},
	\ref{lm:khal6} and \ref{lm:etim} we obtain terms which satisfy the claim.
\end{proof}

\begin{corol}
Let $u\in L$ and $v\in M_{\ge u}$.
	Let $X\in TV\subseteq TV\# H_0$ be a $v$-$[\,]$-word. Then the coproduct $\copr$ of
	$TV\# H_0$ satisfies
	$$
	\copr(X)=X\ot 1 + X_{(-1)}\ot X_{(0)} + \sum_i(X'_ih_i)\ot X''_i,
	$$
	where $X'_i,X''_i\in T^+V$, $h_i\in H_0$, and each $X'_i$ is a $w_i$-$[\,]$-word with
	$w_i\in M_{\ge u}$.
\end{corol}
\begin{proof}
	This follows at once from the Proposition and Lemma~\ref{lm:etim}.
\end{proof}

We study the structure of $R$ now. For that, we will use the following notation.
\begin{defin}
	For $u\in L$, let $\alwu u$ be the subalgebra of $TV$ generated by
	$(\sum_{v\in L,v\ge u}\slwu{v})$. Let $\ilwu u$ be the ideal of $\alwu u$ generated by
	$(\sum_{v\in L,v>u}\slwu{v})$. We define also $\slwl u=\pi(\slwu u)$, $\alwl u=\pi(\alwu u)$,
	$\aalwl u=\pi(\alwu u)\cap \ker \varepsilon $, and
	$\ilwl u=\pi(\ilwu u)$, where $\pi:TV\to R$ is the canonical projection.
\end{defin}
The grading on $R$ induces a grading on all of the algebras and ideals defined above.
Take $u\in L$. Notice that the graded algebras
$\alwu u/\ilwu u$ and $\alwl u/\ilwl u$ have only elements in degrees $m|u|$ where $m\in \ndN _0$.
Moreover, Lemma~\ref{lm:etim} and Corollary~\ref{co:etim2} imply that
$\alwu u$ is the subspace of $TV$ generated by $w$-$[\,]$-words with $w\in M_{\ge u}$
and $\ilwu u$ is the subspace of $\alwu u$ generated by $w$-$[\,]$-words with
$w=w_1\cdots w_n\in M_{\ge u}$ with $w_1>u$.

Notice that the leading vector of any $X\in \alwu u\setminus\ilwu u$ of degree $m|u|$ is a
$u^m$-vector. Thus, we can choose for any $u\in L$ a graded linear map
$\iota _u:\aalwl u/\ilwl u\to T^+\slwu u=\oplus_{m\in\ndN}\slwu {u^m}\subset \alwu u$ such that
$\pi _u\circ \pi |_{T^+\slwu u}\circ \iota _u=\id $, where $\pi :TV\to R$ and
$\pi _u:\alwl u\to \alwl u/\ilwl u$ are the canonical maps.

Recall Def.~\ref{df:op} for the definition of $\ordpr{}$.
\begin{theor}\label{th:isom}
	The map $\phi:\ordpr{u\in L}\aalwl u/\ilwl u\to R^+$ defined by
\stepcounter{equation}
	\begin{equation}
		\label{eq:meq}
		\phi(X_{u_1}\ot \cdots \ot X_{u_n})=
		\pi (\iota_{u_1}(X_{u_1})\cdots\iota_{u_n}(X_{u_n})),
	\end{equation}
	where $u_1,\ldots ,u_n\in L$, $u_1>\cdots >u_n$, and $X_{u_i}\in \aalwl{u_i}/\ilwl{u_i}$ for all $i$,
	is an isomorphism of graded vector spaces.
\end{theor}

To prove the theorem we will use the following lemma.
\begin{lemma}\label{lm:pis}
	Let $u=u_1^{n_1}\cdots u_r^{n_r}\in M$, and let $X_i\in\slwu {u_i^{n_i}}$ for all $i$.
	Take $Y_i=\pi_{u_i}(\pi(X_i))\in\aalwl{u_i}/\ilwl{u_i}$. Then
	$$
	X_1\cdots X_r-\iota_{u_1}(Y_1)\cdots\iota_{u_r}(Y_r) \in \ker\pi
		+ \sum_{\begin{subarray}{c}w>u \\ w\in M\end{subarray}}\slwu w.
	$$
\end{lemma}
\begin{proof}
	Since $\pi_{u_i}\pi(X_i-\iota_{u_i}(Y_i))=0$, we have $X_i-\iota_{u_i}(Y_i)\in\ker\pi+\ilwu{u_i}$.
	We then consider
	\begin{align*}
		Z &= X_1\cdots X_r - \iota_{u_1}(Y_1)\cdots\iota_{u_r}(Y_r) \\
		&= \sum_{j=1}^r X_1\cdots X_{j-1}(X_j-\iota_{u_j}(Y_j))
			\iota_{u_{j+1}}(Y_{j+1})\cdots \iota_{u_{r}}(Y_{r}) \\
		& \in\sum_{j=1}^r\slwu{u_1^{n_1}}\cdots\slwu{u_{j-1}^{n_{j-1}}}
		\ilwu{u_j}\slwu{u_{j+1}^{n_{j+1}}}\cdots\slwu{u_r^{n_r}} + \ker\pi.
	\end{align*}
	As mentioned above the theorem, $\ilwu{u_j}$ consists of sums of
	$w$-$[\,]$-words where $w$ runs on monotonic super-words $>u_j^{n_j}$.  Thus, by
        Corollary~\ref{co:etim2} and Lemma~\ref{lm:dbig},
        $Z$ is a sum of $w$-$[\,]$-words, where $w$ runs on
	super-words $>u$.
\end{proof}

\begin{proof}[Proof of the Theorem.]
	Since the set of words of a given length is finite, the surjectivity of $\phi$ follows easily
	from the previous lemma.

	We now prove injectivity of $\phi $. To do so define $\phi ':\ordpr{u\in L}\aalwl u/\ilwl u
	\to T^+V$ by
	\begin{equation*}
		\label{eq:meq2}
		\phi'(X_{u_1}\ot \cdots \ot X_{u_n})=
		\iota_{u_1}(X_{u_1})\cdots\iota_{u_n}(X_{u_n}),
	\end{equation*}
	where $u_1,\ldots ,u_n\in L$, $u_1>\cdots >u_n$, and $X_{u_i}\in \aalwl{u_i}/\ilwl{u_i}$ for all $i$.
	Assume then that there exists a smallest integer $m$ such that $\phi $ is not injective in degree $m$.
	For all $u\in L$ let $B_u=\{b_{u,i}\,|\,i\in I_u\}$ be a homogeneous basis of $\aalwl u/\ilwl u$,
	where $I_u$ is an appropriate index set, and let $X_{u,i}:=\iota_u (b_{u,i})$ for all $u\in L$,
	$i\in I_u$.

	Suppose that there exists a nonempty finite subset $M'$ of $M$ with $|w|=m$ for $w\in M'$,
	and for each $w=w_1^{n_1}\cdots w_r^{n_r}\in M'$ there exist nonzero elements $b_w\in 
	\aalwl{w_1}/\ilwl{w_1}\ot \cdots \ot \aalwl{w_r}/\ilwl{w_r}$ such that
	$\phi '(\sum _{w\in M'}b_w)\in \ker \pi $.
	Let $u=u_1^{m_1}\cdots u_s^{m_s}$, where $u_1>\cdots >u_s$, be the minimal element of $M'$,
	and write $b_u:=\sum_{i_1,\ldots ,i_s}\lambda _{i_1,\ldots ,i_s}b_{u_1,i_1}\ot \cdots \ot b_{u_s,i_s}$
	with $\lambda _{i_1,\ldots ,i_s}\in \bfi $.
	We consider $TV$ as a subalgebra of $TV\# H_0$, and then we have
	\begin{align*}
	\copr (\phi '(\sum _{w\in M'}b_w)) &=: \sum_iZ'_i\ot Z''_i \\
		&\in (\ker\pi\# H_0)\ot TV+(TV\# H_0)\ot \ker\pi.
	\end{align*}
	Therefore,
	\begin{equation}\label{eq:tsoz}
		\sum_i S^{-1}({Z''_i}_{(-1)})Z'_i\ot {Z''_i}_{(0)} \in \ker\pi\ot TV+TV\ot \ker\pi.
	\end{equation}
	We apply Prop.~\ref{pr:com} to each $[\,]$-letter in $\phi '(\sum _{w\in M'}b_w)$, and we use
	Lemma~\ref{lm:wdec} to obtain a description of the tensor factors of \eqref{eq:tsoz}.
	Afterwards, we apply Corollary~\ref{co:etim2} and Lemma~\ref{lm:dbig} to rearrange the tensor
	factors as sums of $[\,]$-words.
	This gives
	\begin{align*}
		\sum_i &S^{-1}({Z''_i}_{(-1)})Z'_i\ot {Z''_i}_{(0)} \\
			\in& \sum_{i_1,\dots,i_s}\lambda_{i_1,\dots,i_s}
				S^{-1}({X_{u_1,i_1}}_{(-1)}){X_{u_1,i_1}}_{(-2)}X_{u_2,i_2}\cdots X_{u_s,i_s} \ot
				{X_{u_1,i_1}}_{(0)} \\
				&+ \sum_{\begin{subarray}{c}
							w',w''\in M \\
							w'>u_2^{m_2}\cdots u_s^{m_s}\text{ or } w''>u_1^{m_1}
						\end{subarray}}
						\slwu {w'}\ot\slwu {w''}.
	\end{align*}
	By repeatedly using Lemma~\ref{lm:pis} and since $\phi '(\sum _{w\in M'}b_w)\in\ker\pi$, we get
	\begin{align*}
		\sum_i &S^{-1}({Z''_i}_{(-1)})Z'_i\ot {Z''_i}_{(0)} \\
			\in& \sum_{i_1,\dots,i_s}\lambda_{i_1,\dots,i_s}
				\phi'(b_{u_2,i_2})\cdots \phi'(b_{u_s,i_s}) \ot \phi'(b_{u_1,i_1}) \\
				&+ \sum_{\begin{subarray}{c}
							w',w''\in M,\;|w'|,|w''|<|u|, \\
							w'>u_2^{m_2}\cdots u_s^{m_s}\text{ or } w''>u_1^{m_1}
						\end{subarray}}
						(\im\phi'\cap\slwu {w'})\ot(\im\phi'\cap\slwu {w''}) \\
				&+ \ker\pi\ot TV+TV\ot\ker\pi.
	\end{align*}
	Therefore, \eqref{eq:tsoz} shows that
        \begin{equation}\label{eq:kerpi}
	\begin{split}
	\sum_{i_1,\dots,i_s} &\lambda_{i_1,\dots,i_s}
				\phi'(b_{u_2,i_2}\cdots b_{u_s,i_s}) \ot \phi'(b_{u_1,i_1}) \\
				\in& \sum_{\begin{subarray}{c}
							w',w''\in M,\;|w'|,|w''|<|u|, \\
							w'>u_2^{m_2}\cdots u_s^{m_s}\text{ or } w''>u_1^{m_1}
						\end{subarray}}
						(\im\phi'\cap\slwu {w'})\ot(\im\phi'\cap\slwu {w''}) \\
				&+ \ker\pi\ot TV+TV\ot\ker\pi.
	\end{split}
        \end{equation}
	By the assumption on $m$, $\pi\circ\phi'=\phi$ is injective in degrees $<m$ and hence the sums
	\begin{align*}
		& (\ker\pi\cap V^{\ot n})
		+ \Big(\bigoplus_{
				\begin{subarray}{c}
					w'\in M,\;|w'|=n,\\
					w'>u_2^{m_2}\cdots u_s^{m_s}
				\end{subarray}}
			\im\phi'\cap\slwu {w'}\Big)
		+ (\im\phi'\cap\slwu{u_2^{m_2}\cdots u_s^{m_s}} ) \\
		& (\ker\pi\cap V^{\ot n})
		+ \Big(\bigoplus_{
				\begin{subarray}{c}
					w''\in M,\;|w''|=n,\\
					w''>u_1^{m_1}
				\end{subarray}}
			\im\phi'\cap\slwu {w''}\Big)
		+ (\im\phi'\cap\slwu{u_1^{m_1}} )
	\end{align*}
	are direct in $TV$ whenever $1\le n<m$.
	Thus \eqref{eq:kerpi} implies that $\lambda_{i_1,\dots,i_s}=0$ for all $i_1,\dots,i_s$, which contradicts to
	the choice of $b_u$.
\end{proof}

\begin{corol}\label{co:fohs}
	The Hilbert series of $R$ factors as
	\begin{equation*}
		\eta(R,t) = \prod_{u\in L} \eta(\alwl u/\ilwl u, t)
	\end{equation*}
	\hfill\qed
\end{corol}

The importance of Corollary~\ref{co:fohs} becomes clearer with the following theorem.
\begin{theor}\label{th:prna}
	For each $u\in L$, the algebra $\alwl u/\ilwl u\# H_0$ is a $|u|\ndZ$-graded Hopf algebra, where
	the grading is induced by that of $R\# H_0$. Equivalently, $\alwl u/\ilwl u$ is a
	$|u|\ndZ$-graded braided Hopf algebra in $\lYDcat{H_0}$. Moreover, $\alwl u/\ilwl u$ is
	generated by $\slwl u/(\ilwl u\cap \slwl u)$
	and it projects onto the Nichols algebra
	$\cB(\slwl u/(\ilwl u\cap \slwl u))$. The quotient
	$$
		\eta\Big(\alwl u/\ilwl u,t\Big)\;/\;\eta\Big(\cB(\slwl u/(\ilwl u\cap \slwl u)),t^{|u|}\Big)
	$$
	is a power series with nonnegative integer coefficients.
\end{theor}
\begin{proof}
	Since $\alwl u/\ilwl u$ is graded and its degree $0$ part is $H_0$, it is enough to
	prove that $\alwl u/\ilwl u\# H_0$ is a bialgebra (see \cite{a-Takeuchi71}).
	The Theorem follows from Proposition~\ref{pr:whab}, by taking $A=R\# H_0$, $B=\alwl u$, and
	$I=\ilwl u$. It remains to show that Equations~\eqref{eq:abi} hold in this case.
	Indeed, it is enough to prove this for generators of $B$ and $I$, and $[\,]$-letters satisfy
	\eqref{eq:abi} thanks to Proposition~\ref{pr:com}.

	By the definition of $\alwl u$ and $\ilwl u$, $\alwl u/\ilwl u$ is generated as an algebra by
	the space $\slwl u/(\ilwl u\cap \slwl u)$. Further, $\alwl u/\ilwl u$ can be considered as the
	quotient of $T(\slwl u/(\ilwl u\cap \slwl u))$ by a graded Hopf ideal consisting of elements of
	degree $\ge 2$.  Since $\cB(\slwl u/(\ilwl u\cap \slwl u))$ is the quotient of
	$T(\slwl u/(\ilwl u\cap \slwl u))$ by the maximal Hopf ideal consisting of elements of
	degree $\ge 2$, there exists a natural projection $\alwl u/\ilwl u\to \cB(\slwl u/(\ilwl u\cap \slwl u))$.
	The last statement follows from Proposition~\ref{pr:pNa}.
\end{proof}

\begin{remar}
	In Theorem~\ref{th:prna} it is necessary to put $t^{|u|}$ as the variable of the Hilbert series
	of the Nichols algebra, since in $\cB (\slwl u/(\ilwl u\cap \slwl u))$ the elements of
	$\slwl u/(\ilwl u\cap \slwl u)$ are considered to be in degree $1$.
\end{remar}

\begin{oppro}\ 
	\begin{enumerate}
		\item
			Assume that $R$ is a Nichols algebra. Are the graded Hopf algebras $\alwl u/\ilwl u$
			appearing in Theorem~\ref{th:prna} again Nichols algebras? This is true in the case
			where $\mathrm{char}\,\bfi =0$, $H_0$ is the group algebra of an abelian group, and $R$
			is finite dimensional, by Kharchenko's PBW theorem. More generally, if $R$ has a finite
			number of PBW generators, the statement follows by using the Weyl groupoid.
		\item
			Generalize Theorems~\ref{th:isom} and \ref{th:prna} to a more general setting which
			covers also Ufer's PBW basis.
		\item
			Is it possible to generalize results of this paper to arbitrary (say finite dimensional)
			non-semisimple Hopf algebras? 
	\end{enumerate}
\end{oppro}

%\bibliographystyle{amsalpha}
%\bibliography{quantum}

\providecommand{\bysame}{\leavevmode\hbox to3em{\hrulefill}\thinspace}
\providecommand{\MR}{\relax\ifhmode\unskip\space\fi MR }
\providecommand{\MRhref}[2]{%
  \href{http://www.ams.org/mathscinet-getitem?mr=#1}{#2}
}
\providecommand{\href}[2]{#2}

\end{document}